\documentclass{article}

\usepackage{arXiv/PRIMEarxiv}

\usepackage[utf8]{inputenc} 
\usepackage[T1]{fontenc}    
\usepackage{hyperref}       
\usepackage{url}            
\usepackage{booktabs}       
\usepackage{amsfonts}       
\usepackage{amsmath}
\usepackage{nicefrac}       
\usepackage{microtype}      
\usepackage{lipsum}
\usepackage{fancyhdr}       
\usepackage{graphicx}       
\usepackage{siunitx}
\usepackage{algorithm}
\usepackage{algpseudocode}
\usepackage{tikz}

\newcommand{\multifac}{{\it MultiFac}}
\newcommand{\mipseq}{{\it MIP-Seq}}
\newcommand{\lpseq}{{\it LP-Seq}}
\newcommand{\tm}{\textsuperscript{\tt TM}}
\newcommand{\biolog}{Biolog\tm}
\newcommand{\met}[1]{{\tt #1}}
\newcommand{\mysection}[1]{\vspace*{2 mm}{\bf #1}}

\usepackage{natbib}

\graphicspath{{media/}}     

\title{A Combinatorial Optimisation Approach to Multi-factorial
Gap-filling in Genome-scale Metabolic Models (GEMs)
}

\author{
Philip Kilby \\Canberra, Australia, \\ \texttt{pilip@kilby.id.au}
\And
Sevvandi Kandanaarachchi \\ CSIRO Data61, \\ Melbourne, Australia, \\ \texttt{sevvandi.kandanaarachchi@csiro.au}
\And
Matthew J. Morgan \\ CSIRO Environment, \\ Canberra, Australia, \\ \texttt{matthew.morgan@csiro.au}
\And
Amy M. Paten \\ CSIRO Environment, \\ Canberra, Australia, \\ \texttt{amy.paten@csiro.au}
\And
Mariana Velasque  \\ CSIRO Environment, \\ Canberra, Australia, \\ \texttt{mariana.velasque@csiro.au}
\And
Andrew C. Warden \\ CSIRO Environment, \\ Canberra, Australia, \\ \texttt{andrew.warden@csiro.au}
\And
Juan P. Molina Ortiz \\ CSIRO Environment \\ Canberra, Australia, \\ \texttt{jp.molinaortiz@csiro.au}
}

\begin{document}
\maketitle

\begin{abstract}
Genome-Scale Metabolic Models (GEMs) describe the interactions between
genes, proteins, and the biochemical reactions that underpin an
organism's metabolism aiming to computationally simulate functions at
the cellular level. 

While many metabolic reactions can be inferred from genome analysis,
constructing GEMs often involves incorporating reactions unsupported
by genomic data to improve prediction accuracy. This is known as
gap-filling, a process that can be performed manually -- a
time-consuming task -- or computationally.

Traditional computational gap-filling approaches aim to correct GEM
predictions for a single environmental condition (medium) by solving a
large Integer Linear Programming problem. Sequential application
across multiple media can produce a more robust model, but often
introduces unrealistic predictions in other media. They are also slow
to run.

In this paper, we study multi-factorial gap filling, which aims to
gap-fill GEMs across typically 10 or more input media
\emph{simultaneously}, while improving their overall predictive
accuracy and minimising unrealistic behaviour. 

We view the selection of the set of reactions as a combinatorial
optimisation problem, and describe a method based on classic
metaheuristic approaches (Tabu Search and Adaptive Large Neighbourhood
Search) which requires the solution of continuous Linear Programming
problems only.

This paper provides an introduction of this problem to an audience
whose speciality lies outside biology, and suggests a practical
first-cut solution method.

We demonstrate the method on three tests cases, gap-filling GEMs for
the bacteria strains {\it Pseudomonas aeruginosa} PAO1, {\it
Klebsiella pneumoniae} MGH78578  and {\it Escherichia coli} str. K-12
substr. MG1655; selecting 3000 to 4000 reactions from a database of
more than 11000 reactions, while attempting to match the empirically
measured performance on 9 to 28 separate media conditions. We show
that our method outperforms conventional approaches on multiple metrics, 
including Kendal Tau and RMS Error by an average of 7.3\%  and 13.3\%, respectively.
\end{abstract}


\section{Introduction}
An organism's metabolic capabilities are largely determined by its genome, which encodes enzymes that catalyse, sometimes unique, biochemical reactions required to sustain life. Genome-Scale Metabolic Models (GEMs) aim to capture these capabilities as a network of metabolic reactions (\emph{metabolic network}), enabling computational prediction of phenotypes under genetic and nutritional changes, such as growth under defined media. This information can then be applied to strain and bioprocess engineering approaches -- for instance whether a particular change will be effective in promoting more efficient breakdown of plastic \citep{freund2025synthetic}.

A key challenge in constructing a GEM is \emph{gap-filling}: identifying reactions that are used by the organism but are not linked to a known gene (non-gene associated reactions). In the case of microbes, gap-filling is commonly guided by empirical growth measurements across many media conditions. In this paper we use \biolog{} EcoPlate\tm{} growth scores \citep{Garl91Classification}, and optical density (OD) measurements, which provide growth/no-growth information and relative growth strength for each input medium. 

Growth in a defined medium can be simulated using GEMs through Flux Balance Analysis (FBA), which optimises an objective (often biomass production), subject to network capacity constraints. Notably, FBA assumes that the cell is operating at a steady state, meaning that the total amount of produced metabolites equals the total amount of consumed metabolites in terms of mass/stoichiometric balance, and that there is no net accumulation or loss over time. Predicted biomass production can then be compared to empirical growth measurements (EcoPlate, OD) under the same medium.  When predictions do not correlate well with observed growth across one or more media, the reaction set in the current GEM is likely incomplete, motivating gap-filling. Automated gap-filling is often formulated as selecting additional reactions so that growth is feasible in the right medium or better matches observations.

A further consideration is the plausibility that a candidate non-gene associated reaction is used by the organism. This can be informed by evolutionary relatedness: if a close relative is known to use a reaction, it is more likely to be present in the organism under study. We quantify this evidence using \emph{taxonomic distance}, computed from a custom taxonomic tree of a defined set of bacteria (Appendix~\ref{app:tax}), and use it to assign costs to candidate reactions so that reactions supported by close relatives are preferred.

Conventional automated gap-filling approaches typically optimise for feasibility in a single medium (growth vs.\ no-growth) \citep{Hart16Improved}, often applied sequentially across media \citep{agren2013raven}. Sequential gap-filling can be computationally expensive, is order-dependent, and can introduce false positive growth predictions or thermodynamically infeasible cycles (\emph{runaway}). Moreover, beyond basic growth/no growth classification, it usually does not directly target agreement with quantitative growth variation across many media conditions. 

In this paper, we propose a \emph{multi-factorial} gap-filling approach that accounts for multiple media inputs and growth measurements \emph{simultaneously}. We search for a minimum-cost set of reactions, where cost is derived from taxonomic distance, that best matches the empirical performance of an organism across media. Solution quality is evaluated using multiple criteria capturing reaction plausibility and agreement with the measured outcomes, including growth/no-growth consistency and quantitative prediction accuracy.

We treat the problem as a multi-objective discrete optimisation problem and solve it using a metaheuristic search strategy that explores sets of reactions to incorporate into a given GEM. During execution, we seek solutions on the Pareto front induced by the objectives, and the output is a set of candidate GEMs representing different trade-offs between plausibility and predictive agreement. In practice, these outputs are intended to support iterative refinement, allowing domain expertise to be reflected through reaction costs and run parameters.

The main contributions of this paper are:

\begin{itemize}
\item We define the multi-factorial GEM gap-filling problem, and
  provide an introduction to the problem for an audience outside
  biology.
\item We propose objective functions that capture both plausibility of added reactions and agreement with empirical growth outcomes.
\item We present a multi-objective metaheuristic method to solve the resulting discrete optimisation problem.
\item We evaluate this method on three organisms, and compare the
  results with baseline GEMs constructed by a standard method.
\end{itemize}

\section{GEMs and FBA}
\label{sec:GEMs}
A GEM represents an organism's metabolism as a set of biochemical reactions interconnected by metabolites. These models allow us to simulate and study how cells use nutrients and metabolic reactions to grow, respond to genetic or nutritional changes, and produce important compounds (e.g., antibiotics). GEMs are widely used across biotechnology, synthetic biology, and medicine \citep{lu2023silico, molina2022high, marcelino2023disease}. For more details see Appendix~\ref{app:metmod}.

Let $M=\{1,\ldots,m\}$ denote the set of metabolites and $R'=\{0,\ldots,r\}$ the set of reactions, where (without loss of generality) reaction $0$ denotes the biomass production reaction. Each reaction $i \in R'$ is associated with stoichiometric coefficients $a_{ji}$ for metabolites $j \in M$, where $a_{ji}<0$ indicates consumption and $a_{ji}>0$ indicates production. Collecting $a_{ji}$ gives the stoichiometric matrix $A=(a_{ji})$. Reversibility is handled by allowing a reaction to carry flux in either direction (or, equivalently, by splitting into forward and reverse reactions).

Flux Balance Analysis (FBA) predicts feasible reaction fluxes $x_i$ through the metabolic network under a specified medium by enforcing (approximate) steady-state mass balance and reaction capacity constraints while optimising an objective, conventionally \emph{biomass production}. Steady-state mass balance is expressed by $A x \approx 0$ subject to uptake/secretion constraints. In standard formulations, FBA is a linear program \citep{Wats84Metabolic}. For our simulations, we use \emph{parsimonious} FBA (pFBA), which favours solutions achieving growth with minimal overall flux, often used as a proxy for metabolic efficiency \citep{schnitzer2022choice}.

We define $R=\{1,\ldots,r\}$ as the set of reactions excluding the biomass reaction. Let $x_i$ be the flux of reaction $i \in R'$. Medium availability is represented by bounds on net metabolite uptake and secretion. Let $S_j \ge 0$ denote the maximum uptake (supply) for metabolite $j$ and $D_j \ge 0$ denote the maximum secretion or demand for metabolite $j$. Let $U_i$ be an upper bound on the flux of reaction $i \in R$.

We combine biomass maximisation with total-flux minimisation using a scalarisation parameter $\alpha$:
\[
\sum_{i \in R} x_i - \alpha \; x_0.
\]
{\bf pFBA:}
\begin{equation}
\min \sum_{i \in R} x_i - \alpha \; x_0  \label{eq:fbaobj}
\end{equation}
subject to
\begin{align}
-S_j \le \sum_{i \in R'} a_{ji} x_i & \le  D_j \;\;\; \forall j \in
M \label{eq:fba1} \\
x_i & \le U_i  \;\;\; \forall i \in R. \label{eq:fba2}
\end{align}

Constraint (\ref{eq:fba1}) enforces medium uptake and secretion limits while maintaining mass balance, and (\ref{eq:fba2}) bounds reaction fluxes. The parameter $\alpha$ balances biomass production against parsimonious total flux; in practice we increase $\alpha$ until the biomass flux stabilises or a maximum is reached. If biomass remains zero at the maximum $\alpha$, the GEM predicts no growth on the given medium. 

\subsection{Data sources and growth measurements}
\label{sec:data}
In this study, the AGORA2 database \citep{Hein23Genomescale} constitutes the candidate pool of reactions from which gap-filling may select additions. Our empirical targets are previously published growth measurements for three bacterial strains across multiple defined media conditions. Each medium contains a fixed set of components and a designated carbon source; the resulting growth is quantified using \biolog{} EcoPlate\tm{} scores \citep{Garl91Classification} or optical density (OD) measurements. These measurements are approximate proxies for biomass production and are used to define (i) growth vs.\ no-growth and (ii) relative growth strength across media.

Because absolute growth scores are noisy and platform-dependent, our objectives consider both quantitative agreement and the relative ordering of growth across media. In particular, we map measured growth scores to heuristic target biomass production rates by calibrating against the medium with the highest observed growth, as described in Section~\ref{sec:obj}.

During both FBA and conventional culturing, from which OD and Biolog scores are derived, certain metabolites are made available at specific concentrations (culture medium), while others are produced at the end of reaction chains (end-products) or directed to the production of biomass. 

\section{The gap-filling problem}
\label{sec:gapfill}

GEMs assembled from genome annotations can miss reactions due to incomplete or uncertain gene--function associations. \emph{Gap-filling} aims to augment a GEM's reaction set with additional predefined reactions so that simulated growth better agrees with experimental observations across media. This is done by extending FBA to explore which additional reactions better contribute to such objective.

In this work, candidate reactions are drawn from AGORA2 \citep{Hein20AGORA2}. Each candidate reaction $i$ is assigned a cost $RC_i$ reflecting its plausibility for the target organisms, computed from taxonomic distance (Appendix~\ref{app:tax}); lower cost reactions are supported by closer relatives and are preferred.

\subsection{MILP formulation of single media gap-filling problem}
\label{sec:milp}
For a single medium, gap-filling can be posed as selecting reactions to enable/fit growth while penalising implausible additions. We introduce binary variables $y_i$ indicating whether reaction $i \in R$  is available, yielding a mixed-integer formulation:

\[
y_i =
\begin{cases}
1 & \text{if reaction $i$ is enabled} \\
0 & \text{otherwise.}
\end{cases}
\]

{\bf GF:}
\begin{equation}
\min \sum_{i \in R} c_i y_i - \alpha x_0 \label{eq:gfobj}
\end{equation}
subject to
\begin{align}
-S_j \le \sum_{i \in R'} a_{ji} x_i & \le  D_j \;\;\; \forall j \in M \label{eq:gf1} \\
x_i & \le  y_i \, U_i \;\;\; \forall i \in R \label{eq:gf2} \\
y_i & \in  \{ 0, 1 \} \;\;\; \forall i \in R. \label{eq:gfint}
\end{align}

Similar to the pFBA formulation, the objective (\ref{eq:gfobj})
balances the cost of used reactions with the flux through the biomass
equation. Constraint (\ref{eq:gf1}) like constraint (\ref{eq:fba1}),
balances the demand and supply of all metabolites. Constraint
(\ref{eq:gf2}) ensures no flux $x_i$ in reaction $i$ unless $y_i$ is
1, and also enforces the upper bound on flux. Constraint
(\ref{eq:gfint}) ensures $y$ variables are binary.

This formulation is similar to that used by the GlobalFit method
\citep{Hart16Improved}.

\subsection{Cost-weighted pFBA}
\label{sec:fba_c}

When evaluating a candidate reaction set, we use a cost-weighted variant of pFBA to bias flux away from expensive reactions:

{\bf pFBA($c$):}
\begin{equation}
  \min \sum_{i \in R} c_i x_i - \alpha x_0 \label{eq:fbacobj}
\end{equation}
subject to constraints \ref{eq:fba1} and \ref{eq:fba2}.

\section{Multi-factorial gap-filling problem}
\label{sec:mfgf}

Our focus is \emph{multi-factorial} gap-filling: selecting a single reaction set that yields good agreement with empirical growth across \emph{multiple} media simultaneously, while also favouring plausible (low-cost) additions. We evaluate solutions using multiple criteria: reaction cost, growth/no-growth agreement, rank correlation between predicted biomass and growth scores, absolute error (defined formally in Section~\ref{sec:obj}), and mean absolute percentage error (MAPE).

\subsection{Multi-factorial gap-filling formulation}

Given $p$ media indexed by $k \in P=\{1, \ldots,p\}$ with growth scores $G^k$ and medium-specific bounds $S_j^k, D_j^k$, the outer level selects a reaction availability vector $\hat{y} \in \{0,1\}^{|R|}$. For a fixed $\hat{y}$ , we resolve an independent cost-weighted pFBA for each medium:

{\bf pFBA($k,\hat{y},c$):}
\begin{equation}
\min \sum_{i \in R} c_i x^k_i - \alpha^k x^k_0 \label{eq:mffbaobj}
\end{equation}
subject to
\begin{align}
-S^k_j \le \sum_{i \in R'} a_{ji} x^k_i & \le  D^k_j \;\;\; \forall j \in M \label{eq:mffba1} \\
x^k_i & \le  \hat{y_i} \, U_i \;\;\; \forall i \in R. \label{eq:mffba2}
\end{align}

Note that each input medium $k$ has an independent value of the trade-off parameter $\alpha$, $\alpha^k$.
Each medium $k$ has its own trade-off parameter $\alpha^k$. The outer-level search over $\hat{y}$ is solved using the metaheuristic described in Section~\ref{sec:alg}.

\section{Method Description}
\label{sec:method}

\subsection{Objectives}
\label{sec:obj}

For a fixed reaction availability vector $\hat{y}$, we solve pFBA for each input medium $k$ and obtain fluxes $x_i^k$, including biomass flux $x_0^k$. The outer-level objective aggregates four components computed from these solutions 

\mysection{Cost $C(x)$}
\label{sec:cost}

Because the final GEM includes only reactions that carry nonzero flux in at least one medium, we define cost over \emph{used} reactions:
\[
C_u(x) = \sum_{i \in R} RC_i \, \mathbb{I}\left(\max_k x_i^k > 0\right),
\]
where $RC_i$ is the original taxonomic cost (Appendix~\ref{app:tax}). To improve comparability across instances, we scale $C_u(x)$ by a baseline cost $C_0$ obtained by running pFBA with all the optional reactions available:
\[
C(x) = \frac{C_u(x)}{C_0}.
\]

\mysection{Growth Match Error $GME(x)$}
\label{sec:gme}

Let each medium be classified as growth/no-growth using the measurement threshold in Section~\ref{sec:data}, and let predicted growth using the gap-filled GEM be determined by whether $x_0^k$ exceeds a biomass threshold. We define
\[
GME(x) = \left|\left\{k \in P : \text{predicted growth/no-growth does not match observed}\right\}\right|,
\]
so $GME(x) \in \{0,\ldots,p\}$ and lower is better.

\mysection{Correlation $\tau(x)$}

We use Kendall's $\tau$ \citep{Kend38New} to compare the rank ordering of predicted biomass fluxes ($x^k_0$) and measured growth scores ($G^k$). To obtain a minimisation term, we transform it as
\[
\tau'(x) = 1 - \tau(x),
\]
so $\tau'(x)\in[0,2]$.

\mysection{Error $E$}
\label{sec:RMS}

We convert growth scores to heuristic target biomass fluxes by calibrating against a reference medium $l$, preferably inferred from the largest growth score $G^k$, with known biomass flux $x_0^l$ and growth score $G^l$:
\begin{equation}
  \label{eq:targflux}
  \hat{x}_0^k = G^k \cdot \frac{x_0^l}{G^l}.
\end{equation}
The error term $E(x)$ is the root mean square error between $\hat{x}_0^k$ and $x_0^k$ across media

\mysection{Combination and Outer level optimisation}

We combine these objectives using scaling factors $\beta_C, \beta_{GME}, \beta_\tau,$ and $\beta_E$:
\begin{equation}
  \label{eq:outerobj}
  \mathrm{Obj}(x)=\beta_C C(x)+\beta_{GME} GME(x)+\beta_\tau \tau'(x)+\beta_E E(x).
\end{equation}

The outer-level problem can therefore be described as
\begin{equation}
\min_{\hat{y}} \ \beta_C C(x)+\beta_{GME} GME(x)+\beta_\tau \tau'(x)+\beta_E E(x)
\end{equation}
subject to
\begin{equation}
x^k = {\bf pFBA}(k,\hat{y},c) \;\;\; \forall k \in P.
\end{equation}

\subsection{High level overview}
\label{sec:highlev}

The outer level search seeks a reaction set $\hat{y}$ that yields low cost (taxonomically plausible) GEMs whose pFBA biomass predictions ($x_0^k$) agree with measured growth across media under the objectives above. Starting from an initial incumbent $\hat{y}$, each iteration modifies $\hat{y}$ (e.g., enabling/disabling a reaction or updating evaluation costs), resolves pFBA for all media, and updates the incumbent/best solutions and a Pareto archive.

We implement this outer level search using a metaheuristic that combines tabu-style short term memory \citep{Glov93Tabu}, Adaptive Large Neighbourhood Search (ALNS) \citep{Ropk06Adaptive}, and occasional non-improving moves via simulated annealing \citep{Niko10Simulated}. More details about this process are provided in Section~\ref{sec:alg}.

\section{Algorithm description}
\label{sec:alg}

Algorithm~\ref{alg:multifac} summarises \multifac{}. The method maintains a current incumbent GEM represented by reaction availability $\hat{y}\in\{0,1\}^{|R|}$ and (temporary) pFBA costs $c$. Each iteration applies a small modification to $(\hat{y},c)$, resolves pFBA independently for all media, evaluates the multi-component objective, and updates (i) the best solution seen, (ii) the incumbent, and (iii) the Pareto archive.

\begin{algorithm}[htpb]
  \caption{\multifac{} algorithm}\label{alg:multifac}
\begin{algorithmic}[1]
\State Preprocess (set targets/parameters; compute $C_0$)
\State Construct initial incumbent GEM $(\hat{y},c)$
\For{a given number of iterations}
  \State Select an operator (Section~\ref{sec:ops})
  \State Apply operator to update $\hat{y}$ or $c$
  \State Solve \textbf{pFBA}$(k,\hat{y},c)$ for all $k$
  \State Evaluate objective and objective components
  \State Update best-so-far, incumbent acceptance, and Pareto archive
  \State Update tabu restrictions and operator scores; periodically update operator probabilities
\EndFor
\end{algorithmic}
\end{algorithm}

\subsection{Selection algorithms (Operators)}
\label{sec:ops}

The outer-level problem searches over reaction-availability vectors $\hat{y}$ (and, during the run, the evaluation cost vector $c$). We implement this search using an Adaptive Large Neighbourhood Search (ALNS) style framework: a specific \emph{operator} can be selected and applied to iteratively generate a neighbouring candidate solution. In our framework, operators act by selecting a reaction and then excluding it ($\hat{y}_i\!\leftarrow\!0$), re-including it ($\hat{y}_i\!\leftarrow\!1$), or setting a non-unit pFBA cost to unity ($c_i\!\leftarrow\!1$). Tabu restrictions are enforced by excluding recently modified reactions, while gene-indicated reactions are never eligible for exclusion.

When an operator requires choosing a reaction, we use biased roulette-wheel sampling with weights $w_i$ and probability $w_i/\sum_j w_j$ (with $w_i=0$ for ineligible reactions).

Effective operators that have been identified include:

\smallskip\noindent
\textbf{Exclude by cost:} used reactions weighted by $w_i=RC_i$. That is, for reactions $i$ with \mbox{$\max_k x^k_i > 0$}, $w_i =
RC_i$, (0 otherwise).

\noindent
\textbf{Exclude by flux:} used reactions weighted by $w_i=\sum_k x_i^k$ as very high flux reactions are sometimes used in unrealistic pathways/cycles. For all reactions $i$ with \mbox{$\max_k x^k_i > 0$}, $w_i = \sum_k x^k_i$, (0 otherwise).

\noindent
\textbf{Exclude by cost$\times$flux:} used reactions weighted by $w_i=RC_i\sum_k x_i^k$. That is, for reactions $i$ with \mbox{$\max_k x^k_i
> 0$}, $w_i = RC_i \sum_k x^k_i$, (0 otherwise).

\noindent
\textbf{Exclude by overs-vs-unders:} for used reactions, weight equals the count of input media where that
reaction $i$ is used and the input medium is currently \emph{over} target \emph{minus} the count of input media where the reaction is used and the input medium is currently \emph{under} target. If the reaction is used by more ``unders'' than ``overs'', then it is not considered for selection.

\noindent
\textbf{Exclude random:} used reactions weighted uniformly. The weight for each reaction that is currently used is 1.0.

\smallskip\noindent
\textbf{Add by cost:} excluded reactions weighted by $w_i=1/RC_i$. Used to re-introduce a low-cost reaction that has previously been excluded. the weight is the inverse of cost -- i.e. $w_i = \frac{1.0}{RC_i}$. 

\noindent
\textbf{Add random:} excluded reactions weighted uniformly. To explore the solution space without bias, so the weight for each excluded reaction used is 1.0.

\smallskip\noindent
\textbf{Make unit:} if a used reaction has $c_i>1$, select one (weighted by $RC_i$) and set $c_i\leftarrow 1$; if the resulting solution is unacceptable, the reaction is disabled.

\smallskip\noindent
\textbf{Exclude runaway:} if a runaway cycle is detected (i.e., biomass prediction exceeds carbon source supply from media), disable the highest-cost reaction among those attaining the maximum flux.

\subsection{Algorithm components}
\label{sec:components}

\mysection{Preprocessing}

We first run \textbf{pFBA}$(k,\hat{y},c)$ for all $k$ with $\hat{y}\equiv 1$ and $c\equiv RC$. From this run we: (i) record per-medium maximum achievable biomass (an upper bound under subsequent exclusions), updating target biomass values when necessary; (ii) remove media that does not lead to biomass production even with all reactions available; (iii) fix the trade-off parameters $\alpha^k$ for the reminder of the run; and (iv) compute the baseline cost $C_0$ used for cost normalisation (Section~\ref{sec:cost}). 

\mysection{Initial solution}

Starting from gene-indicated reactions only, we consider remaining reactions in increasing $RC_i$ order (ties are broken random). For each candidate reaction $i$, we temporarily enable it and set $c_i\leftarrow 1$, run pFBA for all $k$, and keep $i$ if total biomass does not decrease and no runaway is introduced; otherwise, we disable it. The resulting model is used as the initial incumbent for the improvement loop.

\mysection{Choose a selection algorithm}

Some selection algorithms are enforced when applicable (e.g., \textit{Make unit} when non-unit costs remain; \textit{Exclude runaway} when a runaway is detected). Otherwise, selection algorithms are sampled according to adaptive probabilities (see below). If a selection algorithm cannot be applied (e.g., \textit{Add random} when no reactions are excluded), we resample.

\mysection{Best, incumbent, and Pareto updates}
A new solution replaces the best-so-far if its objective improves. The incumbent is updated if (i) it is new best, (ii) it improves on the incumbent objective, (iii) it has fewer exclusions with similar objective, or (iv) it is accepted under simulated annealing. the Pareto archive is updated when the current GEM dominates an archive member on at least one objective component.

pFBA is run in each turn. Note that if computing resources are available, these can be run in parallel. In our tests we run pFBAs sequentially, but allowed the LP solver to use up to 10 threads simultaneously.

\mysection{Update tabu list}

After a selection algorithm modifies a reaction, reversing that move can be made tabu for a tenure depending on the outcome: use $T_{\text{fail}}$ if $GME$ increases, $T_{\text{worse}}$ if the overall objective worsens, and $T_{\text{incumb}}$ if the solution is accepted as incumbent.

In the first two cases, we do not want to repeat a bad outcome. In the
third case, having discovered an improved GEM, we do not want to
undo it until we have had a chance to explore neighbouring GEMs.

\mysection{Update adaptive selection scores}

Each selection algorithm $k$ accumulates a score $s_k$ and a usage count $n_k$ since the last update: we add $\gamma_1$ for a new-best solution, $\gamma_2$ for incumbent acceptance, and $\gamma_3$ for Pareto-front membership. Every 100 iterations, we update selection probabilities following \citep{Ropk06Adaptive}. For selection algorithm with $n_k>0$ we compute the success measure $m_k=s_k/n_k$ and form $p'_k=m_k/\sum_i m_i$; for $n_k=0$ we keep $p'_k=p_k$. We then apply a learning rate $\sigma$ via $p''_k=\sigma p'_k+(1-\sigma)p_k$ and renormalise $p''$ to obtain the updated probabilities.

\subsection{Implementation notes}

Because the outer loop fixes $\hat{y}$ and only perturbs a small number of values in the stoichiometric matrix between iterations (e.g., toggling $U_i$ for add/exclude or setting $c_i\leftarrow 1$ for make unit), each inner pFBA remains a linear program and can be efficiently re-solved using warm starts; in practice, pFBA sub-problems typically solve in sub-second time.

In our tests we set $\beta_{GME}$ large  so that growth/no-growth agreement dominates the combined objective, and we evaluate multiple $\beta_C,\beta_E$ regimes to explore the cost--error trade-off (Section~\ref{sec:comp}). A workflow diagram of the complete method is provided in Appendix~\ref{app:workflow}

\section{Computational Testing}
\label{sec:comp}

We have named our method \multifac. We performed preliminary
testing to find robust settings for the algorithm's parameters. The
values chosen for testing are given in Table~\ref{tab:params}.

\begin{table}[htbp]
\caption{Algorithm parameters}
  \centering
  \begin{tabular}{l p{10cm} l r}
  Symbol & Meaning & Ref & Value \\
  $T_{\mbox{fail}}$ & Tabu tenure if fewer experiments produce flux &
    Sec~\ref{sec:components} & 200 \\
  $T_{\mbox{worse}}$ & Tabu tenure objective increase &
    Sec~\ref{sec:components} & 25 \\
  $T_{\mbox{incumb}}$ & Tabu tenure if GEM becomes incumbent &
    Sec~\ref{sec:components} & 500 \\
  $\beta_{GME}$ & Multiplier for growth match error component in objective &
    Eq~\ref{eq:outerobj} & 1000 \\
  $\beta_E$ & Multiplier for RMS error component in objective &
    Eq~\ref{eq:outerobj} & Various$^*$ \\
  $\beta_C$ & Multiplier for cost component in objective &
    Eq~\ref{eq:outerobj} &  Various$^*$ \\
  $\beta_\tau$ & Multiplier for $\tau'$, the (negated) correlation component in
    objective & Eq~\ref{eq:outerobj} &  10 \\
  $\gamma_1$ & Reward for new best GEM in adaptive update & 
    Sec~\ref{sec:components}  & 10 \\
  $\gamma_2$ & Reward for new incumbent in adaptive update & 
    Sec~\ref{sec:components}  & 3 \\
  $\gamma_3$ & Reward for Pareto front membership in adaptive update & 
    Sec~\ref{sec:components}  & 1 \\
  $\sigma$ & Learning rate in adaptive update  & 
    Sec~\ref{sec:components}  & 0.3 \\
\end{tabular}
\label{tab:params}
\end{table}

Reaction costs $RC_i$ were calculated using the procedure given in Appendix~\ref{app:tax}

\subsection{Sequential Method}
\label{sec:seq}

To provide reference baselines, we define two sequential gap-filling procedures that process media one-at-a-time (ordered by decreasing growth score). Both baselines share the same outer procedure: for each medium, solve a single medium optimisation problem to select additional reactions; then ``lock in" any newly selected non-gene-indicated reactions by reducing their effective cost (to $1.01$) for subsequent media. The final GEM consists of gene-indicated reactions plus all reactions selected in at least one medium.

\subsection{MIP-based Sequential Method}
\label{sec:mipseq}

\mipseq{} uses the single-medium gap-filling MILP (Section~\ref{sec:milp}) with candidate reactions drawn from AGORA2. Gene-indicated reactions have cost $1.0$ and other reaction costs match those used by \multifac{}. We solve each MILP using the HiGHS solver \citep{Huan18Parallelizing} with 10 threads, a 1-hour time limit per medium and a 1\% MIP-gap termination criterion. Newly selected non-gene-indicated reactions not previously selected are assigned cost $1.01$ for subsequent media (minimum cost for other reactions is $1.5$). 

Although straight-forward, this algorithm has limitations. First, solving the MILP can be quite
time-consuming. We allow one hour per carbon source, and then use the best sub-optimal solution after that time. Additionally, this method is order-dependent and it does not consider biomass flux targets or biomass production ranking. For comparison, we gap-fill carbon sources in order of decreasing \biolog{} score. 

\subsection{LP-based Sequential Method}
\label{sec:lpseq}

\lpseq{} is identical to \mipseq{} except that the per-medium optimisation is based on cost-weighted pFBA LP (Section~\ref{sec:fba_c}) rather than the MILP. This substantially reduces runtime but retains the same limitations of sequential processing. Because both methods maximise biomass for each medium, their biomass outputs (and therefore correlation/RMS terms) are often similar, differing mainly in the reaction sets selected. Due to fluxes being cost-weighted, this method generally selects lower-cost reactions.

\subsection{Data}
\label{sec:data_comp} 

\subsubsection{Test cases}

Three well-studied organisms were selected due to
the availability of a curated GEMs in AGORA2, and publicly available growth scores.

\begin{itemize}
\item {\it Pseudomonas aeruginosa} (strain PA01)
\item {\it Escherichia coli} (strain K12-MG1655)
\item {\it Klebsiella pneumoniae pneumoniae} (strain MGH78578)
\end{itemize}

For each organism, the AGORA2 GEM is used as a baseline comparison. Our goal is to produce gap-filled GEMs that improve over the baseline across the evaluation methods defined in Section~\ref{sec:obj}. Given that AGORA2 GEMs were created using an automated procedure \citep{Hein20AGORA2}, they represent a valid comparison baseline for \multifac{}. 

\subsubsection{Target biomass production rate}
\label{sec:targcalc}

To compute RMS error (Section~\ref{sec:RMS}), we require a target biomass production rate for at least one medium. We take the medium with the highest growth score, run pFBA with all AGORA2 reactions available, and use the resulting biomass flux as the reference. Target biomass rates for other media are then set proportional to their growth scores as in Eq.~\ref{eq:targflux}

\subsection{Run environment}
\label{sec:runenv}

Simulations were run on  Dell PowerEdge C6525 servers with 512 GB of memory. We use the HiGHS linear programming library \citep{Huan18Parallelizing} with 10 threads. \multifac{} is written in C++. Each test was run 10 times for 5000 iterations.

We report results under three objective-weight regimes:
\begin{itemize}
\item \textbf{cost/error:} $\beta_C=100$, $\beta_E=1$ (cost as the primary objective),
\item \textbf{error/cost:} $\beta_E=100$, $\beta_C=1$ (RMS error as the primary objective),
\item \textbf{cost+error:} $\beta_C=1$, $\beta_E=1$ (cost and RMS have equal weights),
\end{itemize}
with $\beta_{GME}=1000$ and $\beta_\tau=10$ fixed as in Table~\ref{tab:params}.

\subsection{Results}

\begin{figure}[htbp]
\centering
\caption{Comparison of methods and performance before and after gap-filling.}
  \fbox{
    \includegraphics[width=1.0\textwidth]{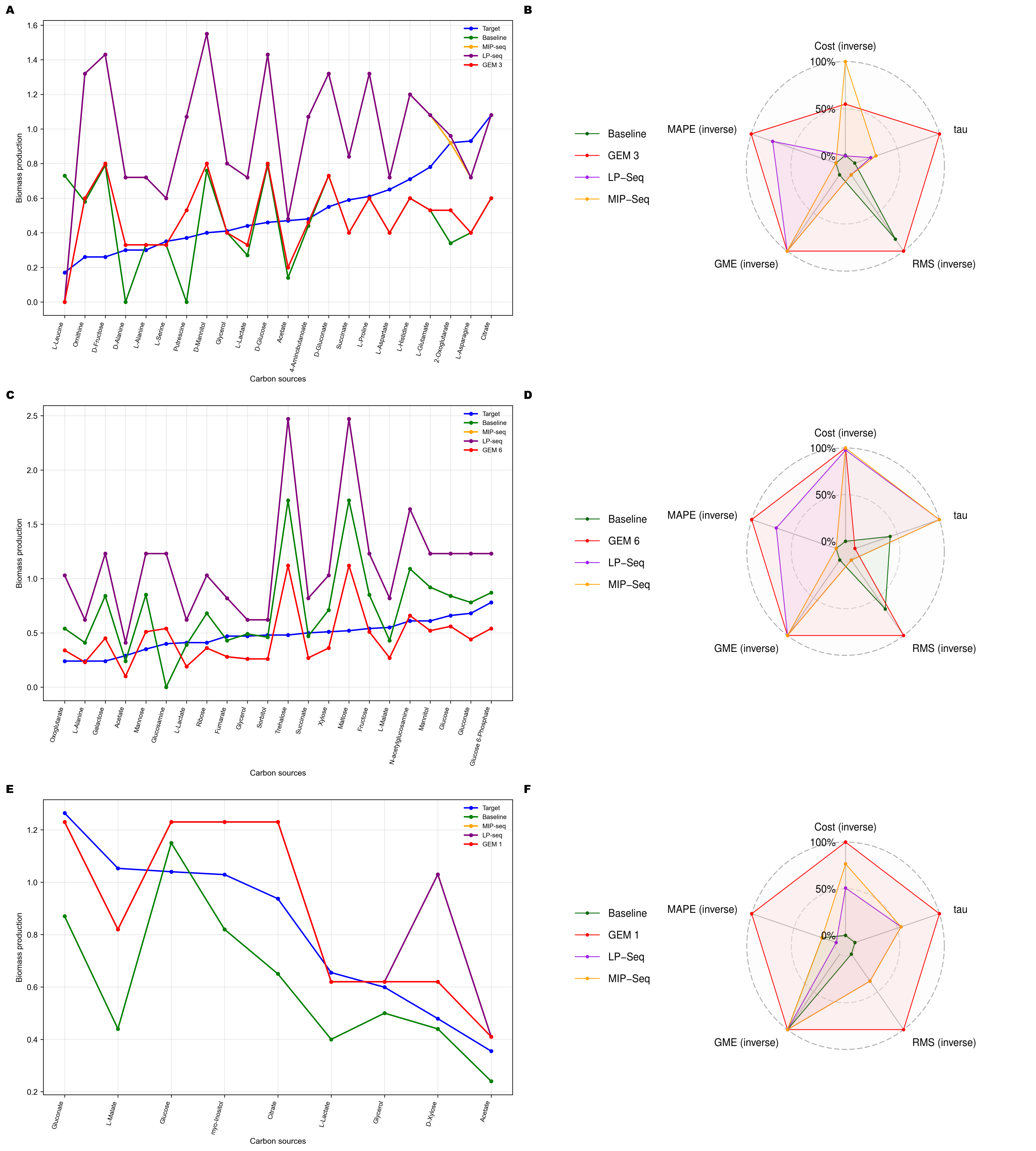}}
    \label{fig:results-panels}
\end{figure}

Figure~\ref{fig:results-panels} summarises prediction-versus-target and multi-metric comparisons for each organism. 

\subsubsection{{\it Pseudomonas aeruginosa}}

We evaluated {\it P. aeruginosa} on 28 \biolog{} EcoPlate carbon sources \citep{oberhardt2008genome}. Six carbon sources were not representable using AGORA2 reactions and could not produce biomass even when all reactions were available; these were excluded, leaving 22 media for analysis. Per-medium targets were obtained by calibrating to the highest-growth medium (citrate), yielding a reference biomass flux of 1.076 (Section~\ref{sec:targcalc}. Detailed per medium targets and predictions are reported in Appendix~\ref{app:growth_pm}.

Table~\ref{tab:pa-results} summarises the Pareto front of GEMs found by \multifac{}, together with the AGORA2 baseline GEM and the two sequential baselines (\mipseq{}, \lpseq{}). The baseline GEM exhibits two false negatives (GME $=2$), whereas all \multifac{} solutions shown achieve perfect growth/no-growth agreement (GME $=0$). Several Pareto-optimal \multifac{} GEMs simultaneously improve upon the baseline in cost, correlation, and RMS error; for example, GEMs~2--4 dominate the baseline across these criteria. As expected, error-focused runs (e.g., GEM~5) achieve the lowest RMS error but at substantially higher reaction cost.

As a representative trade-off we selected GEM~3, which preserves GME $=0$ while improving correlation and RMS error at modest cost compared to the lowest cost solutions. Figure~\ref{fig:results-panels} (panel~A) compares predicted versus target biomass across media for the representative \multifac{} GEM and baselines, and panel~B summarises objective components and Mean Absolute Percentage Error (MAPE), a
measure that \multifac{} was not directly optimising. 

Runtime also favours \multifac{} over MILP-based sequential gap-filling. For {\it P. aeruginosa}, \multifac{} required 20 minutes on average (90\% of runs $<40$ minutes), compared with $>13.5$ hours for \mipseq{} (1--1.9\% MIP gap, average 1.33\%) and approximately 1 minute for \lpseq{}. Consistent with their objective structure, the sequential methods obtain low-cost solutions but perform worse on correlation and RMS error.
\begin{table}[htbp]
\caption{Pareto front of results for the baseline {\it P. aeruginosa} GEM from AGORA2. The
results of the two sequential methods are given in the last three
lines.}
\centering
\begin{tabular}{l S[table-format=6.2] l l l l r }
GEM & Cost & Kendall Tau & RMS Error & GME & Objective   & MAPE \\ \hline
GEM 1 &     670.8  & -0.15  &  0.44  &   0  & cost/error &  \\
GEM 2 &     675.7  &  0.13  &  0.28  &   0  & cost/error &  \\
GEM 3 &     680.1  &  0.25  &  0.28  &   0  & cost+error &  48\% \\
GEM 4 &     686.1  &  0.25  &  0.28  &   0  & cost/error &  \\
GEM 5 &   20498    &  0.61  &  0.23  &   0  & error/cost &  \\ 
GEM 6 &   21159    &  0.65  &  0.23  &   0  & error/cost &  \\
GEM 7 &   23999    &  0.64  &  0.23  &   0  & error/cost &  \\ \hline  
Baseline &  717.6  &  0.09  &  0.33  &   2  &            &  67\% \\
\mipseq  &  649.0  &  0.13  &  0.60  &   0  &            &  123\% \\
\lpseq   &  718.0  &  0.12  &  0.60  &   0  &            &  123\% 
\end{tabular}
\label{tab:pa-results}
\end{table}

\subsubsection{{\it Escherichia coli}}

For {\it E. coli} (K12-MG1655), we used OD-derived growth rates fro 24 carbon sources \citep{adadi2012prediction}. Two media (pyruvate and guanosine) could not produce biomass even with all AGORA2 reactions available and were excluded, leaving 22 media (all classified as growth). Detailed per medium targets and baseline predictions are given in Appendix~\ref{app:growth_pm}.

Table ~\ref{tab:ecoli-results} reports the \multifac{} Pareto front and baseline comparisons. The AGORA2 baseline exhibits one false negative (glucosamine; GME $=1$), whereas all \multifac{} Pareto solutions reported achieve GME $=0$. Compared to the baseline, many \multifac{} solutions reduce RMS error and/or improve rank correlation at substantially lower cost (reflecting avoidance of taxonomically unlikely reactions under out weighting scheme). We selected GEM~6 as a representative compromise between cost and RMS error. Figure~\ref{fig:results-panels} (panels C and D) summarises prediction accuracy and the multi- metric comparison including MAPE. 

In terms of runtime, \multifac{} required 25 minutes on average (90\% $<40$ minutes). \mipseq{} required 20 hours 51 minutes (MIP gap 1--2.8\%, average 1.8\%), while \lpseq{} required approximately 66 seconds. The sequential baselines again show that medium-by-medium optimisation can yield competitive cost and correlation on some instances, but does not consistently minimise RMS error under multi-media targets.
\begin{table}[htbp]
\caption{The full Pareto front for {\it E. coli}
GEMs. The baseline GEM from AGORA2, and the GEMs
created by the sequential methods, form the final three lines.}

\centering
\begin{tabular}{l S[table-format=6.2] l l l l r }
GEM & Cost & Kendall Tau & RMS Error & GME & Objective & MAPE  \\ \hline
GEM 1 &  470.7  &  0.27  &  0.29  &   0  & cost+error  &  \\
GEM 2 &  481.5  &  0.33  &  0.30  &   0  & cost+error  &  \\
GEM 3 &  487.0  &  0.34  &  0.34  &   0  & cost/error  &  \\
GEM 4 &  489.5  &  0.35  &  0.32  &   0  & cost/error  &  \\ 
GEM 5 &  490.4  &  0.34  &  0.28  &   0  & cost/error  &  \\ 
GEM 6 &  493.1  &  0.36  &  0.25  &   0  & cost+error  &  43\% \\ 
GEM 7 &   506.3 &  0.38  &  0.40  &   0  & cost/error  &  \\ 
GEM 8 &   508.5 &  0.39  &  0.73  &   0  & cost/error  &  \\ 
GEM 9 &   21695 &  0.24  &  0.22  &   0  & error/cost  &  \\ 
GEM 10 &  27292 &  0.28  &  0.22  &   0  & error/cost  &  \\ 
GEM 11 & 120178 &  0.48  &  0.17  &   0  & error/cost  &  \\ \hline  
Baseline & 186380 & 0.41 &  0.45  &   1  &             &  75\%  \\
\mipseq  &    487 & 0.48 &  0.82  &   0  &             &  153\% \\
\lpseq   &    552 & 0.48 &  0.82  &   0  &             &  153\% \\
\end{tabular}
\label{tab:ecoli-results} 
\end{table}

\subsubsection{{\it Klebsiella pneumoniae}}

For {\it K. pneumoniae} (MGH78578), we used measured growth rates on 9 carbon sources \citep{liao2011experimentally}; all media were growth and the measured rates were used directly as targets. Table~\ref{tab:kleb-results} summarises the \multifac{} Pareto front and baseline comparisons (per medium values in Appendix~\ref{app:growth_pm}).

Across these media, \multifac{} produces multiple Pareto-optimal GEMs with improved RMS error and/or correlation relative to the AGORA2 baseline. Low-cost solutions are largely represented by the first two Pareto entries, while more expensive solutions enter the front only when they provide a measurable improvement in correlation and/or error. We selected GEM~1 as a representative solution. Figure~\ref{fig:results-panels} (panels~E--F) summarises predicted versus target biomass and the multi-metric comparison including MAPE.

For {\it K. pneumoniae}, \multifac{} required 9 minutes on average (90\% $<18$ minutes), compared with $>9$ hours for \mipseq{} (MIP gap 1.7--3.2\%, average 2.5\%) and approximately 27 seconds for \lpseq{}. As in other cases, sequential baselines tend to prioritise low-cost feasibility but do not consistently improve correlation and RMS error under multi-media evaluation.

\begin{table}[htbp]
\caption{The Pareto front of GEMs for {\it Klebsiella pneumoniae}. The pre-gap-filling AGORA2 GEM is the
last line.}

  \centering
  \begin{tabular}{l S[table-format=6.2] l l l l r }
GEM & Cost & Kendall Tau & RMS Error & GME & Objective & MAPE \\ \hline
GEM 1 &    375.0  &  0.72  &  0.15  &   0  & cost+error  & 16\% \\
GEM 2 &    448.6  &  0.72  &  0.14  &   0  & cost+error  & \\
GEM 3 &    14993  &  0.61  &  0.12  &   0  & error/cost  & \\
GEM 4 &    58835  &  0.50  &  0.12  &   0  & error/cost  & \\
GEM 5 &    65960  &  0.61  &  0.12  &   0  & error/cost  & \\ 
GEM 6 &   106026  &  0.67  &  0.09  &   0  & error/cost  & \\ \hline  
Baseline & 642.2  &  0.61  &  0.29  &   0  &             & 28\% \\
\mipseq  & 437.5  &  0.67  &  0.24  &   0  &             & 26\% \\
\lpseq   & 506.5  &  0.67  &  0.24  &   0  &             & 26\% 
\end{tabular}

\label{tab:kleb-results}

\end{table}

\section{Conclusions}
\label{sec:conc}

Genome-scale metabolic models (GEMs) are valuable predictive tools, but their utility depends on whether they reproduce observed growth across diverse media conditions. This paper formulated \emph{multi-factorial} gap-filling: selecting a single set of reactions that matches growth targets across multiple carbon sources \emph{simultaneously}, rather than medium-by-medium. Considering all media jointly avoids key limitations of sequential gap-filling, including order dependence and the tendency to prioritise feasibility over agreement  with quantitative growth variation.

We introduced objective components that capture complementary aspects of solution quality, including reaction plausibility (via taxonomically derived costs), growth/no-growth agreement, rank consistency, and quantitative error. We then presented \multifac{}, a practical heuristic that searches over reaction availability decisions while solving only linear pFBA subproblems per iteration. The method is multi-objective and returns a Pareto set of candidate GEMs, allowing users to choose models that best trade off plausibility and predictive fit for their application.

Across three organisms, \multifac{} produced GEMs that matched or improved upon the AGORA2 baseline and two sequential baselines on the evaluation measures, while remaining computationally tractable on standard hardware. In particular, the simultaneous multi media formulation enabled models with perfect growth/no-growth agreement and improved quantitative fit in cases where baseline and sequential approaches showed false negatives or higher prediction error. 

Future work includes evaluating alternative plausibility models (e.g., phylogenetic similarity) and integrating complementary constraints (e.g., enzyme and proteome allocation dynamics) to further reduce unrealistic flux distributions and improve predictive accuracy.

In practice, we view \multifac{} as a tool for \emph{iterative} model development: the Pareto set provides draft GEMs that can be refined by adjusting reaction costs, restricting candidate reaction pools, or incorporating additional biological evidence. Overall, \multifac{} provides a practical and scalable route to generating high-quality draft GEMs that better reflect multi-condition growth data prioritising plausible biology.

\bibliographystyle{unsrt}  
\bibliography{multifac}

\appendix
\section{Technical details}

\subsection{Calculation of taxonomy-based reaction costs}
\label{app:tax}

To guide gap-filling toward more biologically plausible models,
\multifac{} uses a taxonomy-informed weighting scheme that integrates
reaction frequency information derived from AGORA2 GEMs across
hierarchical groupings of genomes. The frequency of every candidate
reaction is computed from the reference GEMs at progressively broader
taxonomic levels (species, genus, family, order, class, phylum) from
the genome of the target strain. This provides a quantitative measure
of how frequently a reaction is observed in increasingly distant
relatives of the target strain. A monotonic mapping is then applied to
convert these frequencies into costs, such that reactions that occur
frequently in close relatives receive low weights, while those
observed only rarely or only in distant groups are assigned higher
weights.

For each reaction, the minimum cost across all taxonomic levels is
retained and incorporated into the gap-filling optimisation as a
reaction-specific penalty. This procedure biases models towards adding
reactions that have been observed in taxonomically similar
organisms. This reduces the inclusion of biologically implausible
reactions and produces more robust, evolutionarily grounded GEMs.

The algorithm  description is given in Algorithm~\ref{alg:tax}.

The algorithm has as input a database of reference GEMs that includes
\begin{itemize}
\item the taxonomic classification of each organism
\item which reactions are included in the model for the organism
\end{itemize}
We then also require the taxonomic classification of the target organism.

Such data is available, for example, from the GTDB-tk database
\citep{Parks22GTDB}.

In biology, taxonomic classification is a hierarchical scheme divided into
\emph{ranks}. From most specific classification to broadest, the
hierarchy of ranks is
\begin{itemize}
\item Species,
\item Genus, 
\item Family, 
\item Order, 
\item Class
\item Phylum,
\item Kingdom,
\item Domain
\end{itemize}

For example, the Red Fox is 
\begin{itemize}
\item species vulpes,
\item genus Vulpes,
\item family Canidae,
\item order Carnivora,
\item class Mammalia,
\item phylum Chordata,
\item kingdom Animalia,
\item domain Eukaria.
\end{itemize}

In the weighting scheme used here, ranks \emph{kingdom} and
\emph{domain} are not considered, as AGORA2 only encompasses bacteria and a few archaea.

If an organism from the database and the target organism match in
\emph{species}, then they are very close relatives. If they match
first in \emph{family} (that is, their species and genus are
different, but their family is the same), then they are more distant
relatives -- for instance your pet Labrador matches the \emph{family}
of the Red Fox (\emph{Canidae}), but is from genus \emph{Canis}. If they
match in first in \emph{phylum}, then there are many millions of years
of evolution separating the two organisms.

The weighting method relies on rank-specific cost ranges. If the
target matches the \emph{species}, then the species cost range is
used. If the target matches first in \emph{genus} then the genus cost
range is used, etc.

\begin{tabular}{l r r}
  Rank & Lower cost $l_k$ & Upper Cost $u_k$ \\
species &  1  	& 2  	\\ 
genus 	&  2  	& 10  	\\ 
family 	&  10  	& 50  	\\ 
order 	&  50  	& 250 	\\ 
class 	&  250  & 1250  \\ 
phylum 	&  1250 & 6250  
\end{tabular}

In the algorithm below, 
\begin{itemize}
\item organism $T$ is the target organism
\item We have a predicate \emph{matches}$(m, T, k)$ which is true if
  the taxonomy of the organism for which $m$ is a model matches the
  target organism $T$ at rank $k$.
\item $m^k$ is the number of models that match the target at rank $k$.
\item $p^k_i$ is proportion of models that use reaction $i$ in rank $k$.
\item $w_i$ is the reaction cost used for reaction $i$.
\end{itemize}

For example, given a model \emph{lab} for a Labrador, \\
\emph{matches} (\emph{lab}, Red Fox, \emph{family}) is {\tt true}, but \\
\emph{matches} (\emph{lab}, Red Fox, \emph{genus}) is {\tt false}

\begin{algorithm}[htpb]
  \caption{The taxonomic weight calculation algorithm}\label{alg:tax}
  \begin{algorithmic}[1]
    \State{{\bf Input:} $T$, the target organism; $DB$ the database of
    taxonomy-indexed GEMs}
    \For{rank $k$ {\bf in} (species, genus, family, order, class,
      phylum)}
      \State{$S^k \gets \{m \in DB \mid$ \emph{matches}$(m, T, k)\}$}
      \State{$m^k \gets \mid S^k \mid$}
      \For{each reaction $i \in R$}
        \Comment{Find the proportion of models that use this reaction}
        \State{$n^k_i \gets$ the number of times reaction $i$ is used by a
          model in $S^k$}
        \State{$p^k_i = \frac{n^k_i}{m^k}$}
      \EndFor
    \EndFor

    \For{each reaction $i \in R$}
      \For{rank $k$ {\bf in} (species, genus, family, order, class, phylum)}
        \State{} \Comment{Use an inverse linear scale to calculate the cost for
          this rank according to the proportion $p^k_i$}
        \State{$w^k_i \gets l_k + (1 - p^k_i) * (u_k - l_k)$}
      \EndFor
      \State{$w_i = \min_k{w^k_i}$}
    \EndFor
      
    \State{{\bf return:} The vector $(w_i : i \in R)$}
  \end{algorithmic}
\end{algorithm}

\section{Workflow diagram}
\label{app:workflow}

\begin{figure}[htbp]
  \caption{The complete \multifac{} algorithm.}
    \fbox{\includegraphics[width=0.95\textwidth]{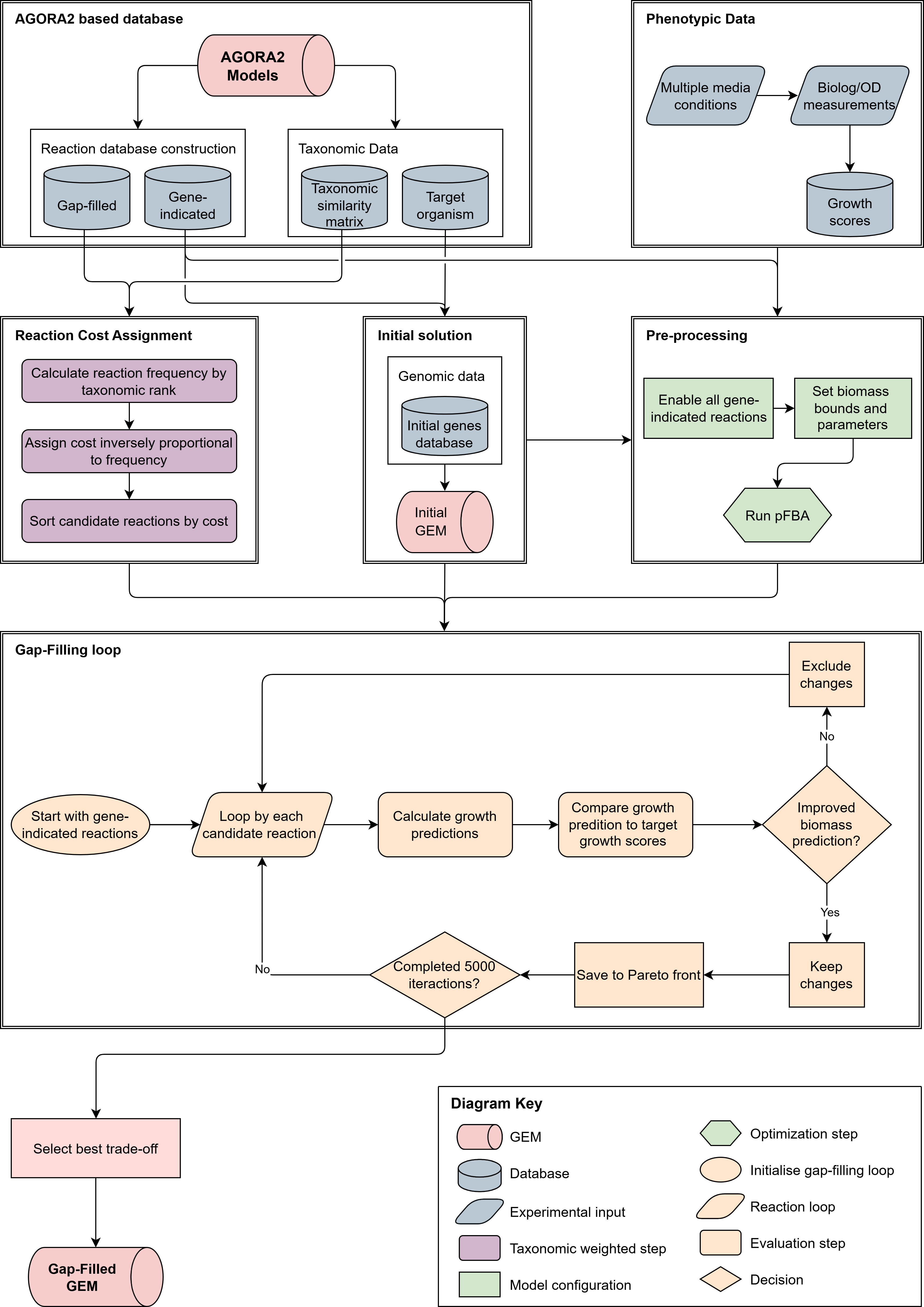}}
  \label{fig:workflow}
\end{figure}

For readers that prefer a visual representation of the algorithm, we
include a workflow diagram as Figure~\ref{fig:workflow}.

\section{Metabolic modelling basics}
\label{app:metmod}

Metabolic modelling is the theoretical and computational framework that uses GEMs, among other tools, to study metabolism. GEMs can be utilised to perform metabolic modelling through the application of constraints (e.g., limiting the concentration of nutrients available) and optimisation algorithms such as FBA. The result of such analyses is reported as \emph{flux distributions}, where \emph{flux} is the rate at which metabolites flow through each metabolic reaction, predicting how a cell distributes resources to achieve a goal. Conventionally flux in flux distributions is expressed in units of millimoles of metabolite per gram of cell dry weight per hour (mmol/gDW/h). 

\subsection{Elements of a GEM}
\label{app:basic}

The basic elements of a GEM are \emph{metabolites} and
\emph{reactions}, both of which are conventionally assigned unique identifiers. In this study, we utilise reaction and metabolite identifiers from AGORA2 \citep{Hein23Genomescale}, a database of over 7000 GEMs and about 11000 reactions. In this context, metabolites are the building blocks
of reactions. For example, the reaction ``Aminobutyraldehyde
Dehydrogenase'' involves the metabolites described on Table~\ref{tab:metabolites}

\begin{table}[ht]
\caption{Metabolites participating in the reaction ``Aminobutyraldehyde
Dehydrogenase'' and their chemical formulas}
\centering
\begin{tabular}{l l l}
\hline
\textbf{Metabolite ID} & \textbf{Name} & \textbf{Chemical Formula} \\
\hline
\met{abut}  & 4-Aminobutanoate & $\mbox{C}_4 \mbox{H}_9 \mbox{N} \mbox{O}_2$ \\
\met{abutn} & 4-Aminobutanal & $\mbox{C}_4 \mbox{H}_{10} \mbox{NO}$ \\
\met{h2o}   & Water & $\mbox{H}_2 \mbox{O}$ \\
\met{nad}   & Nicotinamide adenine dinucleotide &
$\mbox{C}_{21} \mbox{H}_{26} \mbox{N}_7 \mbox{O}_{14} \mbox{P}_2$ \\
\met{h}     & Hydrogen & H \\
\met{nadh}  & Reduced nicotinamide adenine dinucleotide &
$\mbox{C}_{21} \mbox{H}_{27} \mbox{N}_7 \mbox{O}_{14} \mbox{P}_2$ \\
\hline
\end{tabular}
\label{tab:metabolites}
{}
\end{table}

The reaction is written as
\[
\mbox{\met{abutn}} + \mbox{\met{h2o}} + \mbox{\met{nad}} \rightarrow
\mbox{\met{abut}} + 2 \; \mbox{\met{h}} + \mbox{\met{nadh}}
\]

The amount of each metabolite involved in the reaction is called the
\emph{stoichiometric coefficient}. In the above equation, the
stoichiometric coefficients are all 1 (unit coefficients are not
usually shown) except for hydrogen (\met{h}), which has a coefficient
of 2.  Thus, for a unit flux through this reaction, 1 unit each of
\met{abutn}, \met{h2o} and \met{nad} is consumed, while one unit of
\met{abut}, one unit of \met{nadh} and two units of \met{h} are
produced. Metabolites that are consumed in a reaction are called
``reactants'', and those that are produced are called ``products''.

In the AGORA2 database, reactants are represented with negative
stoichiometric coefficients and products with positive coefficients,
so the above reaction would appear as \\
-1 \met{abutn} -1 \met{h2o} -1 \met{nad} 1 \met{abut} 2 \met{h} 1 \met{nadh}

Reactions in GEMs often are reversible; meaning they can run in the forward or reverse direction. If the above reaction
were marked as reversible, then the reaction
\[
\mbox{\met{abut}} + 2 \; \mbox{\met{h}} + \mbox{\met{nadh}} \rightarrow
\mbox{\met{abutn}} + \mbox{\met{h2o}} + \mbox{\met{nad}}
\]

would also be available, and the reaction notation ``$\leftrightarrow$''
would be used.

\subsection{Metabolic network}
\label{app:metnet}  

In GEMs reactions that use a common metabolite can be linked to form a
network, as shown in Figure~\ref{fig:metnet}. The reactions labelled
\textbf{R1}, \textbf{R2} and \textbf{R3} take the input metabolites
\textbf{A} and \textbf{B}, and convert them into \textbf{F} and
\textbf{G}, producing and consuming intermediate metabolites \textbf{C}, \textbf{D}
and \textbf{E} in the process. Individual stoichiometric reactions are
arranged into a stoichiometric matrix, with negative coefficients
indicating reaction reactants and positive coefficients indicating
reaction products. This is done to guarantee mass balance based on metabolite chemical formulae.

\begin{figure}[htbp]
  \fbox{
    \includegraphics[width=1.0\textwidth]{images/metnet.png}
  }
  \caption{A simple metabolic network, using
      reactions R1, R2 and R3, and metabolites A through G. Stoichiometric
      coefficients for each reaction are organised in a mass balanced
      stoichiometric matrix.}
  \label{fig:metnet}
\end{figure}

\section{Measured and predicted growth rate comparison tables}
\label{app:growth_pm}

\begin{table}[htbp]
\caption{Data for {\it P. aeruginosa}}
  \centering
  \begin{tabular}{l l l l l l l }
Carbon & \biolog & Growth & Target  & Baseline  & \mipseq  & \lpseq \\
Source & score   & class  & biomass & biomass   &  biomass &  biomass \\ \hline
L-Leucine &         143	 & No-growth &   0.17  &  0.73          &   0.00  &  0.00  \\
Ornithine &         219	 & Growth    &   0.26  &  0.58          &   1.32  &  1.32  \\
D-Fructose &        221	 & Growth    &   0.26  &  0.79          &   1.43  &  1.43  \\
D-Alanine &         249	 & Growth    &   0.30  &  \textbf{0.00} &   0.72  &  0.72  \\
L-Alanine &         252	 & Growth    &   0.30  &  0.33          &   0.72  &  0.72  \\
L-Serine &          298	 & Growth    &   0.35  &  0.33          &   0.60  &  0.60  \\
Putrescine &        313	 & Growth    &   0.37  &  \textbf{0.00} &   1.07  &  1.07  \\
D-Mannitol &        337	 & Growth    &   0.40  &  0.76          &   1.55  &  1.55  \\
Glycerol &          342	 & Growth    &   0.41  &  0.40          &   0.80  &  0.80  \\
L-Lactate &         370	 & Growth    &   0.44  &  0.27          &   0.72  &  0.72  \\
D-Glucose &         385	 & Growth    &   0.46  &  0.79          &   1.43  &  1.43  \\
Acetate &           393	 & Growth    &   0.47  &  0.14          &   0.48  &  0.48  \\
4-Aminobutanoate &  403	 & Growth    &   0.48  &  0.44          &   1.07  &  1.07  \\
D-Gluconate &       464	 & Growth    &   0.55  &  0.73          &   1.32  &  1.32  \\
Succinate &         495	 & Growth    &   0.59  &  0.40          &   0.84  &  0.84  \\
L-Proline &         510	 & Growth    &   0.61  &  0.60          &   1.32  &  1.32  \\
L-Aspartate &       547	 & Growth    &   0.65  &  0.40          &   0.72  &  0.72  \\
L-Histidine &       597	 & Growth    &   0.71  &  0.60          &   1.20  &  1.20  \\
L-Glutamate &       652	 & Growth    &   0.78  &  0.53          &   1.08  &  1.08  \\
2-Oxoglutarate &    774	 & Growth    &   0.92  &  0.34          &   0.92  &  0.96  \\
L-Asparagine &      778	 & Growth    &   0.93  &  0.40          &   0.72  &  0.72  \\
Citrate &           905	 & Growth    &   1.08  &  0.60          &   1.08  &  1.08  
\end{tabular}
\label{tab:pa01-cs}
\end{table}

\begin{table}[htbp]
\caption{Data for {\it E. coli}}
\centering
\begin{tabular}{l l l l l l }
Carbon & Growth & Target  & Ref. GEM & \mipseq & \lpseq  \\
Source & class  & biomass & biomass    & biomass & biomass \\ \hline
Oxoglutarate            & Growth &   0.24  &   0.54  &        	1.03 & 1.03 \\
L-Alanine               & Growth &   0.24  &   0.41  &        	0.62 & 0.62 \\
Galactose               & Growth &   0.24  &   0.84  &        	1.23 & 1.23 \\
Acetate                 & Growth &   0.29  &   0.24  &        	0.41 & 0.41 \\
Mannose                 & Growth &   0.35  &   0.85  &        	1.23 & 1.23 \\
Glucosamine             & Growth &   0.40  &   \textbf{0.00} &	1.23 & 1.23 \\
L-Lactate               & Growth &   0.41  &   0.39  &        	0.62 & 0.62 \\
Ribose                  & Growth &   0.41  &   0.68  &        	1.03 & 1.03 \\
Fumarate                & Growth &   0.47  &   0.43  &        	0.82 & 0.82 \\
Glycerol                & Growth &   0.47  &   0.49  &        	0.62 & 0.62 \\
Sorbitol                & Growth &   0.48  &   0.46  &        	0.62 & 0.62 \\
Trehalose               & Growth &   0.48  &   1.72  &        	2.47 & 2.47 \\
Succinate               & Growth &   0.50  &   0.47  &        	0.82 & 0.82 \\
Xylose                  & Growth &   0.51  &   0.71  &        	1.03 & 1.03 \\
Maltose                 & Growth &   0.52  &   1.72  &        	2.47 & 2.47 \\
Fructose                & Growth &   0.54  &   0.85  &        	1.23 & 1.23 \\
L-Malate                & Growth &   0.55  &   0.43  &        	0.82 & 0.82 \\
N-acetylglucosamine     & Growth &   0.61  &   1.09  &        	1.64 & 1.64 \\
Mannitol                & Growth &   0.61  &   0.92  &        	1.23 & 1.23 \\
Glucose                 & Growth &   0.66  &   0.84  &        	1.23 & 1.23 \\
Gluconate               & Growth &   0.68  &   0.78  &        	1.23 & 1.23 \\
Glucose 6-Phosphate     & Growth &   0.78  &   0.87  &        	1.23 & 1.23
\end{tabular}
\label{tab:ecoli-cs}
\end{table}

\begin{table}[htbp]

\caption{Data for {\it Klebsiella pneumoniae}}
\centering
\begin{tabular}{l l l l l l }
  Carbon & Growth & Target  & Ref. GEM & \mipseq & \lpseq  \\
Source & class  & biomass & biomass      & biomass & biomass \\ \hline
Gluconate	 & Growth & 1.264 &   0.87 	 &   1.23 &  1.23 \\
L-Malate	 & Growth & 1.053 &   0.44 	 &   0.82 &  0.82 \\
Glucose		 & Growth & 1.04  &   1.15 	 &   1.23 &  1.23 \\
myo-Inositol & Growth & 1.029 &   0.82   &   1.23 &  1.23 \\
Citrate		 & Growth & 0.937 &   0.65   &   1.23 &  1.23 \\
L-Lactate	 & Growth & 0.655 &   0.40 	 &   0.62 &  0.62 \\
Glycerol	 & Growth & 0.599 &   0.50 	 &   0.62 &  0.62 \\
D-Xylose	 & Growth & 0.479 &   0.44 	 &   1.03 &  1.03 \\
Acetate		 & Growth & 0.355 &   0.24 	 &   0.41 &  0.41 
\end{tabular} 
\label{tab:kleb-cs}
\end{table}

\end{document}